\documentclass[a4paper,11pt,english]{smfart}

\usepackage{amssymb}
\usepackage{url}
\usepackage[T1]{fontenc}
\usepackage[cal=boondoxo]{mathalfa}

\usepackage{lscape}
\usepackage{fancybox}
\usepackage{todonotes}

\usepackage[leftbars]{changebar}

\usepackage{framed}

\usepackage{palatino}
\usepackage{rotating}
\usepackage{graphicx}

\usepackage{enumerate}

\usepackage{color}
\definecolor{shadecolor}{gray}{0.90}

\def\bfit{\bfseries\itshape}

\def\proofname{D\'emonstration}

\input xypic
\xyoption{all}
\xyoption{arc}

\makeindex

\def\indexname{Index des notations}

\newtheorem{theo}{Theorem}[section]
\def\thetheo{\thesection.\arabic{theo}}
\newtheorem{prop}[theo]{Proposition}

\newtheorem{lem}[theo]{Lemma}

\newtheorem{defi}[theo]{Definition}

\newtheorem{conj}[theo]{Conjecture}

\renewcommand\theequation{\thetheo}
\def\equat{\refstepcounter{theo}\begin{equation}}
\def\endequat{\end{equation}}

\renewcommand\thetable{\thetheo}

\renewcommand\thesection{\arabic{section}}
\renewcommand\thesubsection{\thesection.\Alph{subsection}}

\def\contentsname{Table des mati\`eres}
\def\listfigurename{Liste des figures}
\def\listtablename{Liste des tables}
\def\appendixname{Appendice}

\def\refname{R\'ef\'erences}

  \def\bG{{\mathfrak b}}  
\def\CG{{\mathfrak C}}    \def\CM{{\mathbb{C}}}

  \def\gG{{\mathfrak g}}

  \def\lG{{\mathfrak l}}

    \def\RM{{\mathbb{R}}}

    \def\ZM{{\mathbb{Z}}}

    \def\CC{{\mathcal{C}}}
    \def\DC{{\mathcal{D}}}

\def\Ob{{\mathbf O}}

\def\Crm{{\mathrm{C}}}

\def\Mrm{{\mathrm{M}}}

\def\Zrm{{\mathrm{Z}}}

\def\a{\alpha}

\def\g{\gamma}
\def\G{\Gamma}

\def\D{\Delta}
\def\e{\varepsilon}

\def\l{\lambda}
\def\L{\Lambda}

\def\o{\omega}
\def\O{\Omega}
\def\r{\rho}
\def\s{\sigma}

\def\th{\theta}

\def\t{\tau}
\def\x{\xi}
\def\z{\zeta}

\def\mub{{\boldsymbol{\mu}}}

\DeclareMathOperator{\End}{{\mathrm{End}}}

\DeclareMathOperator{\Id}{{\mathrm{Id}}}

\DeclareMathOperator{\Ind}{{\mathrm{Ind}}}

\DeclareMathOperator{\Irr}{{\mathrm{Irr}}}

\DeclareMathOperator{\Reg}{{\mathrm{Reg}}}

\DeclareMathOperator{\Tr}{{\mathrm{Tr}}}

\def\to{\rightarrow}
\def\longto{\longrightarrow}

\def\DS{\displaystyle}

\def\SSS{\scriptscriptstyle}

\def\finl{~$\blacksquare$}

\def\lexp#1#2{\kern\scriptspace\vphantom{#2}^{#1}\kern-\scriptspace#2}
\def\le{\hspace{0.1em}\mathop{\leqslant}\nolimits\hspace{0.1em}}
\def\ge{\hspace{0.1em}\mathop{\geqslant}\nolimits\hspace{0.1em}}

\mathchardef\inferieur="321E
\mathchardef\superieur="321F

\def\eqna{\begin{eqnarray*}}
\def\endeqna{\end{eqnarray*}}

\def\maxi{{\mathrm{max}}}

\def\mini{{\mathrm{min}}}
\def\itemth#1{\item[${\mathrm{(#1)}}$]}

\catcode`\@=11
\long\def\@car#1#2\@nil{#1}
\long\def\@first#1#2{#1}
\long\def\@second#1#2{#2}
\long\def\ifempty#1{\expandafter\ifx\@car#1@\@nil @\@empty
  \expandafter\@first\else\expandafter\@second\fi}
\catcode`\@=12

\DeclareMathOperator{\REF}{Ref}

\def\boitegrise#1#2{\begin{centerline}{\fcolorbox{black}{shadecolor}{~
    \begin{minipage}[t]{#2}{\vphantom{~}#1\vphantom{$A_{\DS{A_A}}$}}
            \end{minipage}~}}\end{centerline}\medskip}

\theoremstyle{remark}
\newtheorem{rema}[theo]{Remark}

\newtheorem{exemple}[theo]{Example}

\theoremstyle{plain}

\def\BIL{LR}
\def\GAUCHE{L}
\def\CAR{CAR}
\def\FAM{FAM}

\def\reg{{\mathrm{reg}}}

\def\calo{\SSS{\Crm\Mrm}}

\def\xyinj{\ar@{^{(}->}}
\def\xysur{\ar@{->>}}

\bigskip\def\unb{{\boldsymbol{1}}}

\def\kl{\SSS{\mathrm{KL}}}

\makeatletter
\def\hlinewd#1{%
\noalign{\ifnum0=`}\fi\hrule \@height #1 %
\futurelet\reserved@a\@xhline}
\makeatother

\newlength\epaisLigne

\usepackage{multirow}

\makeatletter
\def\hlinewd#1{%
\noalign{\ifnum0=`}\fi\hrule \@height #1 %
\futurelet\reserved@a\@xhline}
\makeatother

\usepackage{multirow}

\usepackage{lscape}

\usepackage{pstricks}

\def\spectrum{{\mathcal{S\!p}}}

\begin{document}


\title{Calogero-Moser cells of dihedral groups at equal parameters}

\author{{\sc C\'edric Bonnaf\'e}}
\address{IMAG, Universit\'e de Montpellier, CNRS, Montpellier, France} 
\email{cedric.bonnafe@umontpellier.fr}

\makeatother

\author{{\sc J\'er\^ome Germoni}}

\address{Universit\'e de Lyon, Universit\'e Claude Bernard Lyon 1, CNRS UMR 5208,
Institut Camille Jordan, F-69622 Villeurbanne, France}
\email{germoni@math.univ-lyon1.fr}

\date{\today}

\thanks{The first author is partly supported by the ANR:
Projects No ANR-16-CE40-0010-01 (GeRepMod) and ANR-18-CE40-0024-02 (CATORE)}

\begin{abstract}
We prove that Calogero-Moser cells coincide with Kazhdan-Lusztig cells 
for dihedral groups in the equal parameter case.
\end{abstract}

\maketitle
\pagestyle{myheadings}

\markboth{\sc C. Bonnaf\'e \& J. Germoni}{\sc Calogero-Moser cells of dihedral groups}

\bigskip

Calogero-Moser cells have been defined by Rouquier and the first author 
for any finite complex reflection group and any parameter, based on ramification theory for Calogero-Moser 
spaces~\cite{pacific, calogero}. It is conjectured that, for Coxeter groups, 
Calogero-Moser cells coincide with Kazhdan-Lusztig 
cells~\cite[Conj.~3.1~and~3.2]{pacific},~\cite[Conj.~LR~and~L]{calogero}, which were 
defined by Kazhdan-Lusztig~\cite{KL} in the equal parameter case and by Lusztig~\cite{lusztig} 
in the general case. 
The aim of this paper is to prove this conjecture for dihedral groups 
in the equal parameter case.

For Calogero-Moser {\it left} cells, an alternative (and partially conjectural) 
definition is proposed in~\cite[Theo.~13.3.2]{calogero}, based on Gaudin operators. 
This definition is recalled in Section~\ref{sec:gaudin}. 
This is the point of view we adopt in this paper: in the relatively small case 
of dihedral groups, 
an explicit diagonalization of these operators is possible, and the 
computation of Calogero-Moser left cells becomes easy. 


\bigskip

\section{Setup}

\medskip

Let $V$ be a finite dimensional Euclidean real vector space, whose 
positive definite symmetric bilinear form is denoted by $(.,.)$ and let 
$W$ be a finite subgroup of the orthogonal group $\Ob(V)$ generated 
by reflections. For $v \in V$, we denote by $v^*$ the element of the dual space $V^*$ 
defined by $v^*(y)=(y,v)$ for all $y \in V$. The map $V \to V^*$, $v \mapsto v^*$ 
is a $W$-equivariant isomorphism of vector spaces. 

The set of reflections of $W$ is denoted by $\REF(W)$. 
For $\a \in V \setminus \{0\}$, we denote by $s_\a$ the orthogonal 
reflection such that $s_\a(\a)=-\a$. We set
$$\Phi=\{\a \in V~|~(\a,\a)=1\text{~and~}s_\a \in W\}.$$
Then $\Phi=-\Phi$, and we fix a subset $\D$ of $\Phi$ of cardinality $\dim \RM\Phi$ 
such that every element of $\Phi$ belongs to $\sum_{\a \in \D} \RM_{\geqslant 0} \,\a$ 
or to $\sum_{\a \in \D} \RM_{\leqslant 0} \,\a$. We set
$$S=\{s_\a~|~\a \in \D\},$$
so that $(W,S)$ is a finite Coxeter system. We set 
$$\Phi^+= \Phi \cap \sum_{\a \in \D} \RM_{\geqslant 0} \,\a\qquad\text{and}\qquad \Phi^-=-\Phi^+,$$
so that $\Phi=\Phi^+ ~\dot{\cup}~ \Phi^-$, where $\dot{\cup}$ means disjoint union. 
We set 
$$v_0=\sum_{\a \in \Phi^+} \a.$$
We denote by $w_0$ the longest element of $W$ (with respect to the length 
function $\ell : W \to \ZM_{\geqslant 0}$ defined by the choice of $S$). Then
\equat\label{eq:w0}
w_0(v_0)=-v_0.
\endequat
We set
$$V_\reg=V \setminus \bigcup_{\a \in \Phi} V^{s_\a}$$
$$\CG=\{v \in V~|~\forall~\a \in \D,~(\a,v) > 0\}.\leqno{\text{and}}$$
Then $\CG$ is the {\it fundamental chamber} of $W$ associated with $S$ and $v_0 \in \CG$. 
Recall that its closure is a fundamental domain for the action of $W$ on $V$.

We denote by $\Reg_W$ the character afforded by the regular representation and 
$\Irr(W)$ denotes the set of irreducible characters of $W$. We denote 
by $\unb_W$ the trivial character of $W$ and we set $\e : W \to \mub_2$, $w \mapsto \det(w)$. 
We denote by $\CC_\RM$ the vector space of maps 
$c : \REF(W) \longto \RM$ such that $c_s=c_t$ if $s$ and $t$ are conjugate in $W$ 
(the elements of $\CC_\RM$ are called {\it parameters}). 
Finally, if $X$ is a subset of $W$, we set $X^{-1}=\{w^{-1}~|~w \in W\}$. 

\bigskip

\section{Recollection about Kazhdan-Lusztig cells}

\medskip

Let $c \in \CC_\RM$. 
To the datum $(W,S,c)$ 
are associated three partitions of $W$ into {\it Kazhdan-Lusztig left}, {\it right} and 
{\it two-sided $c$-cells} (see for instance~\cite[Chap.~6]{livre}). 
To each Kazhdan-Lusztig left $c$-cell $C$ is 
associated a {\it Kazhdan-Lusztig $c$-cellular character} that is denoted by $\chi_C^{c,\kl}$. 
Then
\equat\label{eq:reg-kl}
\Reg_W =\sum_{C} \chi_C^{c,\kl},
\endequat
where $C$ runs over the set of Kazhdan-Lusztig left $c$-cells.

On the other hand, to each Kazhdan-Lusztig two-sided $c$-cell $\G$ of $W$ is associated a subset 
$\Irr_\G^{c,\kl}(W)$ called the {\it Kazhdan-Lusztig $c$-family} associated with $\G$. 
They form a partition of $\Irr(W)$:
\equat\label{eq:two}
\Irr(W) =\dot{\bigcup_\G} ~\Irr_\G^{c,\kl}(W),
\endequat
where $\G$ runs over the set of Kazhdan-Lusztig two-sided $c$-cells. 
Here are some other properties of Kazhdan-Lusztig cells 
(see for instance~\cite[\S{6.1},~\S{6.2}~and~Chap.~10]{livre}).

\bigskip

\begin{prop}\label{prop:kl-w0}
Let $C$ (resp. $\G$) be a Kazhdan-Lusztig left (resp. two-sided) $c$-cell. Then:
\begin{itemize}
\itemth{a} $C^{-1}$ is a Kazhdan-Lusztig right $c$-cell and $|C|=\chi_C^{c,\kl}(1)$.

\itemth{b} $\G$ is a union of Kazhdan-Lusztig left (or right) 
$c$-cells. Moreover,
$$|\G|=\sum_{\chi \in \Irr_\G^{c,\kl}(W)} \chi(1)^2.$$

\itemth{c} If $C \subset \G$, then every irreducible 
component of $\chi_C^{c,\kl}$ belongs to $\Irr_\G^{c,\kl}(W)$.

\itemth{d} $w_0\G w_0=\G$.

\itemth{e} $Cw_0$ and $w_0C$ (resp. $w_0\G=\G w_0$) are Kazhdan-Lusztig left (resp. two-sided) 
$c$-cells. Moreover,
$$\chi_{C w_0}^{c,\kl}=\chi_{w_0 C}^{c,\kl} = \chi_C^{c,\kl} \cdot \e
\qquad\text{and}\qquad \Irr_{\G w_0}^{c,\kl}(W)=\Irr_\G^{c,\kl}(W)\cdot \e.$$

\itemth{f} If $c_s \neq 0$ for all $s \in \REF(W)$, then $\{1\}$ and $\{w_0\}$ 
are Kazhdan-Lusztig two-sided $c$-cells. If moreover $c_s > 0$ for all $s \in \REF(W)$, then 
$$
\begin{cases}
\chi_{\{1\}}^{c,\kl}=\e,\\ 
\chi_{\{w_0\}}^{c,\kl}=\unb_W,
\end{cases}\qquad\text{and}\qquad
\begin{cases}
\Irr_{\{1\}}^{c,\kl}(W)=\{\e\},\\
\Irr_{\{w_0\}}^{c,\kl}(W)=\{\unb_W\}.
\end{cases}
$$

\itemth{g} If $\t : W \to \mub_2$ is a linear character, then $C$ (resp. $\G$) 
is a Kazhdan-Lusztig left (resp. two-sided) $\t \cdot c$-cell. Morever 
$$\chi_C^{\t \cdot c,\kl}=\chi_C^{c,\kl} \cdot \t\qquad\text{and}\qquad
\Irr_\G^{\t \cdot c,\kl}(W) = \Irr_\G^{c,\kl}(W) \cdot \t.$$
\end{itemize}
\end{prop}

\bigskip

In the above statement~(g), 
$\t\cdot c$ denotes the element of $\CC_\RM$ defined by $(\t\cdot c)_s=\t(s)c_s$.

\bigskip

\section{Gaudin operators}

\medskip

For $y \in V$, $v \in V_\reg$ and $v' \in V$, we define an endomorphism 
$D_y^{c,v,v'}$ of the underlying vector space of 
the group algebra $\RM W$ by the following formula~\cite[\S{13.2}]{calogero}:
$$\forall~w \in W,\quad D_y^{c,v,v'}(w)=(y,w^{-1}(v')) w -\sum_{\a \in \Phi^+} c_{s_\a} 
\cfrac{(y,\a)}{(v,\a)} ws_\a.$$
The endomorphism $D_y^{c,v,v'}$ is called a {\it Gaudin operator} 
(and is somewhat similar to Dunkl operators (see for instance~\cite[3.1.B]{calogero}). 
Then the map $D^{c,v,v'} : V \to \End_\RM(\RM W)$ is linear and 
it follows from~\cite[\S{13.2}]{calogero}\footnote{Note that we have not used exactly the convention 
of~\cite[\S{13.2}]{calogero}: our operators are obtained from those in loc.~cit. by 
conjugating by the $\RM$-linear map extending the inversion $w \mapsto w^{-1}$ in $W$ and 
by identifying $V$ and $V^*$ thanks to the non-degenerate form $(.,.)$.} 
that 
\equat\label{eq:commute}
\bigl[ D_y^{c,v,v'},D_{y'}^{c,v,v'} \bigr]=0
\endequat
for all $y$, $y' \in V$, $v \in V_\reg$ and $v' \in V$. Now, for $\l \in V^*$, we set 
$$E_\l^{c,v,v'}=\{e \in \RM W~|~D_y^{c,v,v'}(e)=\l(y) e\}$$
and we define
$$\spectrum^{c,v,v'}=\{\l \in V^*~|~E_\l^{c,v,v'} \neq 0\}.$$
As all reflections of $W$ have order $2$, the matrix of $D_y^{c,v,v'}$ in the canonical 
basis of $\RM W$ is real and symmetric, so it is diagonalizable. Therefore, for all $(v,v') \in V_\reg \times V$, 
the family of commuting 
matrices $\DC^{c,v,v'}=(D_y^{c,v,v'})_{y \in V}$ is simultaneously diagonalizable. In other words,
\equat\label{eq:decomposition}
\RM W =\bigoplus_{\l \in \spectrum^{c,v,v'}} E_\l^{c,v,v'}
\endequat
for any $(v,v') \in V_\reg \times V$. 
The set $\spectrum^{c,v,v'}$ is called the {\it spectrum} of the family $\DC^{c,v,v'}$. 
We say that the family $\DC^{c,v,v'}$ has {\it simple spectrum} if $|\spectrum^{c,v,v'}|=|W|$ 
(in other words, if $\dim E_\l^{c,v,v'}=0$ or $1$ for all $\l \in V^*$). 

\bigskip

\begin{quotation}
\begin{conj}\label{conj:simple}
If $c \in \CC_\RM$ and $(v,v') \in V_\reg \times V_\reg$, then the family $\DC^{c,v,v'}$ has simple spectrum.
\end{conj}
\end{quotation}

\bigskip

By the work of Mukhin-Tarasov-Varchenko~\cite[Coro.~7.4]{MTV1},~\cite{MTV2}, 
this conjecture is known to hold in type $A$. Here is a weaker 
form of this conjecture. 

\bigskip

\begin{quotation}
\begin{conj}\label{conj:simple-facile}
If $c \in \CC_\RM$ and $\xi$, $\xi' \in \RM_{>0}$, 
then the family $\DC^{c,\xi v_0,\xi' v_0}$ has simple spectrum.
\end{conj}
\end{quotation}

\bigskip

We will prove in this paper that this weaker form holds if $W$ is dihedral and $c$ is constant 
(which is the so-called ``equal parameter case'').

\bigskip

\begin{exemple}\label{ex:t=0}
The matrix of the endomorphism $D_y^{0,v,v'}$ in the canonical basis 
of $\RM W$ is diagonal, and so its spectrum can be easily computed. 
We get
$$\spectrum^{0,v,v'} = \{w(v^{\prime *})~|~w \in W\}.$$
In particular, $\DC^{0,v,v'}$ has simple spectrum if and only if $v' \in V_\reg$.\finl
\end{exemple}

\bigskip

We conclude this subsection by some relations between Gaudin operators. 
For $w \in W$, we denote by $l_w$ (resp. $r_w$) the automorphism of the 
$\RM$-vector space $\RM W$ defined by left (resp. right) 
multiplication by $w$ (resp. $w^{-1}$). If $\t : W \to \mub_2$ is a linear character, 
we denote by $\t_\bullet$ the automorphism of the $\RM$-algebra $\RM W$ 
defined by $\t_\bullet(w)=\t(w)w$ for all $w \in W$. The following 
formulas are straightforward:
\equat\label{eq:left right linear}
\begin{cases}
l_w D_y^{c,v,v'} l_w^{-1}=D_y^{c,v,w(v')},\\
\\
r_w D_y^{c,v,v'} r_w^{-1}=D_{w(y)}^{c,w(v),v'},\\
\\
\t_\bullet D_y^{c,v,v'} \t_\bullet^{-1}=D_y^{\t\cdot c,v,v'}.
\end{cases}
\endequat

\bigskip

\section{Calogero-Moser cells}\label{sec:gaudin}

\medskip

\subsection{Calogero-Moser cellular characters}
The operator $D_y^{c,v_0,0}$ commutes with left multiplication by 
$\RM W$. So each subspace $E_\l^{c,v_0,0}$ inherits a structure 
of $\RM W$-module: we denote by $\chi_\l^c$ the character afforded by this 
$\RM W$-module. We define the {\it Calogero-Moser $c$-cellular characters} 
to be the characters of the form $\chi_\l^c$ for some 
$\l \in \spectrum^{c,v_0,0}$. Note that we may have $\chi_\l^c=\chi_\mu^c$ 
even if $\l \neq \mu$. Then~\eqref{eq:decomposition} implies that 
\equat\label{eq:reg}
\Reg_W=\sum_{\l \in \spectrum^{c,v_0,0}} \chi_\l^c.
\endequat
In particular, every irreducible character of $W$ occurs in some 
Calogero-Moser $c$-cellular character.

Replacing $(c,v_0)$ by $(\xi c, \xi' v_0)$ (with $\xi$, $\xi' \in \RM^\times$) 
amounts to multiplying the Gaudin operators by $\xi/\xi'$: this does not change the 
list of Calogero-Moser cellular characters. This shows that 
Calogero-Moser $c$-cellular characters coincide with 
Calogero-Moser $\xi c$-cellular characters.

\bigskip

\begin{rema}\label{rem:pas simple}
The family $\DC^{c,v_0,0}$ does not have a simple spectrum in general. Indeed, 
if $W$ is not abelian, then an irreducible character of degree $> 1$ 
occurs in some cellular character $\chi_\l^c$, which shows that 
$\dim E_\l^{c,v_0,0} \ge 2$.\finl
\end{rema}

\bigskip

\subsection{Left cells}
In order to define Calogero-Moser left cells, we need to work under the 
following hypothesis:

\medskip

\boitegrise{{\bf Hypothesis.} {\it In this subsection, and only in this subsection, 
we assume that Conjecture~\ref{conj:simple-facile} holds. }}{0.75\textwidth}

\medskip

\noindent Let $v_1$, $v_2 \in \RM_{>0} v_0$. We fix two continuous functions
$\g$, $\xi : [0,1] \longto \RM_{\geqslant 0}$ 
such that $\g(t) \ge 0$ and $\xi(t) > 0$ for all $t \in [0,1)$ and
$$
\begin{cases}
\g(0)=0,\quad \xi(0)=1,\\
\g(1)=1,\quad \xi(1)=0.
\end{cases}
$$
Therefore, for $t \in [0,1)$, the family $\DC^{\g(t)c,v_1,\xi(t)v_2}$ has simple spectrum 
(indeed, if $\g(t)=0$, then this follows from Example~\ref{ex:t=0} and, if 
$\g(t) > 0$, then $D_y^{\g(t)c,v_1,\xi(t)v_2}=D_y^{c,\g(t)^{-1} v_1,\xi(t)v_2}$ 
and so this follows from the fact that we assume that Conjecture~\ref{conj:simple-facile} holds). 
So this spectrum varies continuously according to the parameter $t$. 
But, for $t=0$, we have $\spectrum^{0,v_1,v_2} = \{w(v_2^*)~|~w \in W\}$ by Example~\ref{ex:t=0}. 
This means that, for each $w \in W$, there exists a unique continuous 
map $\l_w : [0,1] \to V^*$ such that 
$$
\begin{cases}
\l_w(0)=w(v_2^*)\\
\text{$\l_w(t) \in \spectrum^{\g(t)c,v_1,\xi(t)v_2}$ for all $t \in [0,1]$.}\\
\end{cases}
$$
and the family $(\l_w)_{w \in W}$ satisfies that 
\equat\label{eq:chemains disjoints}
\forall~t \in [0,1),~\l_w(t) \neq \l_{w'}(t)
\endequat
whenever $w \neq w'$. 
However, it may happen that $\l_w(1) = \l_{w'}(1)$ even if $w \neq w'$. This leads 
to the following definition:

\bigskip

\begin{defi}\label{defi:cm-cells}
Two elements $w$ and $w'$ are said to belong to the same 
{\bfit Calogero-Moser left $c$-cell} if $\l_w(1)=\l_{w'}(1)$.

If $C$ is a Calogero-Moser left $c$-cell, we set $\chi_C^{c,\calo}=\chi_{\l_w(1)}^c$ 
(where $w$ is some, or any, element of $C$): it is called 
the {\bfit Calogero-Moser $c$-cellular character} associated with $C$.
\end{defi}

\bigskip

\begin{rema}\label{rem:souplesse}
A simple choice would be to take $v_1=v_2=v_0$, $\g(t)=tc$ and $\xi(t)=1-t$. 
But we want to work with this slightly more general setting for 
more flexibility. Indeed, one could wonder whether the notion of Calogero-Moser 
left $c$-cell depends on the choices of $v_1$, $v_2$, $\g$, $\xi$. In fact, 
it does not, because the topological space $\CC_\RM \times \RM_{>0} \times \RM_{>0}$ 
is simply connected.

For instance, this shows that, if $r \in \RM_{>0}$, then Calogero-Moser left $rc$-cells 
coincide with Calogero-Moser left $c$-cells, and their associated cellular 
characters agree.

If we assume moreover that Conjecture~\ref{conj:simple} holds, then we could have 
added some more flexibility, by taking $v_1$, $v_2$ in $\CG \times \CG$ and 
replaced the path $t \mapsto \xi(t) v_2$ by any path $\nu_2 : [0,1] \to \overline{\CG}$ 
such that $\nu_2(t) \in \CG$ for $t \in [0,1)$, $\nu_2(0)=v_2$ and $\nu_2(1)=0$ and the path 
$t \mapsto \g(t)c$ by any path $[0,1] \to \CC_\RM$ starting at $0$ and ending at $c$.\finl
\end{rema}

\bigskip

The Formula~\eqref{eq:reg} can be rewritten as follows:
\equat\label{eq:reg-calo}
\Reg_W=\sum_{C} \chi_C^{c,\calo},
\endequat
where $C$ runs over the set of Calogero-Moser left $c$-cells.

The following conjecture has been proposed in~\cite[Conj.~3.2]{pacific} and~\cite[Conj.~L]{calogero}:

\bigskip

\begin{quotation}
\begin{conj}\label{conj:kl=cm}
Calogero-Moser left $c$-cells coincide with Kazhdan-Lusztig left $c$-cells. Moreover, 
if $C$ is one of these, then $\chi_C^{c,\calo}=\chi_C^{c,\kl}$.
\end{conj}
\end{quotation}

\bigskip

A very weak evidence for this Conjecture is the comparison between~\eqref{eq:reg-kl} 
and~\eqref{eq:reg-calo}. Note also that it holds for $c=0$, as easily 
shown in~\cite[Coro.~17.2.3]{calogero}. 
A somewhat strong evidence for this Conjecture is that it holds in type $A$, 
by the work of Brochier-Gordon-White~\cite{BGW}. The aim of this paper is to deal with 
the far easier (but still non-trivial) case of dihedral groups whenever $c$ is constant. The 
following list of 
properties of Calogero-Moser left cells shows that Conjecture~\ref{conj:kl=cm} is compatible with 
Proposition~\ref{prop:kl-w0}. 

\bigskip

\begin{prop}\label{prop:calo-w0}
Let $C$ be a Calogero-Moser left $c$-cell. Then:
\begin{itemize}
\itemth{a} $|C|=\chi_C^{c,\calo}(1)$.

\itemth{b} $Cw_0$ and $w_0C$ are Calogero-Moser left 
$c$-cells. Moreover,
$$\chi_{C w_0}^{c,\calo}=\chi_{w_0 C}^{c,\calo} = \chi_C^{c,\calo} \cdot \e.$$

\itemth{c} If $c_s \neq 0$ for all $s \in \REF(W)$, then $\{1\}$ and $\{w_0\}$ 
are Calogero-Moser left $c$-cells. If moreover $c_s > 0$ for all $s \in \REF(W)$, then 
$$
\chi_{\{1\}}^{c,\calo}=\e\qquad\text{and}\qquad 
\chi_{\{w_0\}}^{c,\calo}=\unb_W.
$$

\itemth{d} If $\t : W \to \mub_2$ is a linear character, then $C$ (resp. $\G$) 
is a Calogero-Moser left (resp. two-sided) $\t \cdot c$-cell. Morever 
$$\chi_C^{\t\cdot c,\calo}=\chi_C^{c,\calo} \cdot \t.$$
\end{itemize}
\end{prop}

\bigskip

\begin{proof}
As explained in Remark~\ref{rem:souplesse}, we may assume that $v_1=v_2=v_0$ and 
that $\g(t)=tc$ and $\xi(t)=1-t$ for all $t \in [0,1]$.

\medskip

(a) is clear.

\medskip

(b) Let $\t_0 = \e_\bullet \circ l_{w_0} : \RM W \longto \RM W$. 
Since $\e\cdot c=-c$ and $w_0(v_0)=-v_0$, we get 
from~\eqref{eq:left right linear} that
$$\t_0 D_y^{tc,v_0,(1-t)v_0}\t_0^{-1}=D_y^{-tc,v_0,(t-1)v_0}=D_{-y}^{tc,v_0,(1-t)v_0}.$$
This means that $\l \in \spectrum^{tc,v_0,(1-t)v_0}$ if and only if 
$-\l  \in \spectrum^{tc,v_0,(1-t)v_0}$. Since $\l_{ww_0}(v_0)=ww_0(v_0)=-w(v_0)=-\l_w(0)$, this 
shows that $\l_{ww_0}(t)=-\l_w(t)$ for all $t \in [0,1]$. In particular, $Cw_0$ 
is a Calogero-Moser left $c$-cell.

Finally, if $\l \in \spectrum^{c,v_0,0}$, then 
$E_{-\l}^{c,v_0,0} = \t_0(E_\l^{c,v_0,0})=\e_\bullet(E_\l^{c,v_0,0})$. 
This proves that $\chi_{-\l}^c = \chi_\l^c \cdot \e$, and completes the proof of~(a).

\medskip

(d) follows from the third equality in~\eqref{eq:left right linear} and the same 
argument as in~(b). 

\medskip

(c) By using~(d) and rectifying the signs if necessary thanks to a linear 
character, we may, and we will, assume that $c_s > 0$ for all $s \in \REF(W)$. 
We have
$$D_{-v_0}^{tc,v_0,(1-t)v_0}(w)= -(t-1) ( v_0,w^{-1}(v_0) ) w + \sum_{\a \in \Phi^+} tc_{s_\a} ws_\a.$$
Let $A$ denote the diagonal endomorphism $D_{-v_0}^{0,v_0,v_0}$ and let $B$ denote the 
Gaudin operator $D_{-v_0}^{c,v_0,0}$, so that $D_{-v_0}^{tc,v_0,(1-t)v_0} = (1-t)A + t B$. 

The matrix $B$ is a real matrix with non-negative coefficients, which is 
primitive (because $W$ is generated by $\REF(W)$). Let $\nu=(v_0,v_0)+1$. 
Then $A + \nu \Id_{\RM W}$ is a diagonal matrix with positive coefficients. Therefore, 
if $t > 0$, the matrix $D_{-v_0}^{tc,v_0,(1-t)v_0} + \nu \Id_V$ is a real matrix 
with non-negative coefficients which is primitive. By Perron-Frobenius Theorem, 
its spectral radius $\r_t$ is an eigenvalue of $D_{-v_0}^{tc,v_0,(1-t)v_0} + \nu \Id_V$, 
with multiplicity $1$. Therefore, $\r_t$ varies continuously as $t$ varies. 

For $t=0$, $A + \nu \Id_{\RM W}$ is diagonal and its biggest diagonal coefficient 
is $(v_0,v_0) + \nu$, which occurs with multiplicity $1$ (and its eigenvector is $w_0$). 
So the map $\r' : [0,1] \to \RM_{>0}$, $t \mapsto \r_t-\nu$ is continuous. 
Adding $\nu \Id_{\RM W}$ was an artefact to obtain a matrix with non-negative 
coefficients and apply Perron-Frobenius Theorem. Coming back to 
$D_{-v_0}^{tc,v_0,(1-t)v_0}$, we have proven that 
$\r_t'$ is its biggest eigenvalue, that it has multiplicity $1$, 
that it varies continuously for $t \in [0,1]$, that $\r_0'=(v_0,v_0)$ 
with eigenvector $w_0$.

Let us now prove that $\r_1'=a$, where $a=\sum_{\a \in \Phi^+} c_{s_\a}$. For this, note 
first that $\sum_{w \in W} w$ is an eigenvector of $B$ with eigenvalue $a$: this 
proves that $\r_1' \ge a$. Now, 
if $\o=\sum_{w \in W} p_w w$ is an eigenvector of $B$ for the eigenvalue $\r_1'$, 
then Perron-Frobenius Theorem says that $p_w > 0$ for all $w \in W$. 
Let $w_1 \in W$ be such that $p_{w_1}$ is maximal. Then the coefficient 
of $w_1$ in $B(\o)$ is equal to $\sum_{\a \in \Phi^+} c_{s_\a} p_{w_1s_\a} \le a p_{w_1}$. 
But this coefficient if $\r_1' p_{w_1}$. This proves that $\r_1' \le a$. 
Consequently, $\r_1'=a$ and the corresponding eigenvector is $\sum_{w \in W} w$. 

This shows that $\l_{w_0}(t)(-v_0)=\r_t'$ and so $\l_{w_0}(1)$ has multiplicity $1$ 
in $\spectrum^{0,v_0,v_0}$ and the corresponding eigenspace $E_{\l_{w_0}(1)}^{0,v_0,v_0}$ 
is the line generated by $\sum_{w \in W} w$. 
This concludes the proof of the fact that $\{w_0\}$ is alone in its 
Calogero-Moser left $c$-cell, and that the associated cellular character is $\unb_W$. 
The statement for $\{1\}$ instead of $\{w_0\}$ is now obtained by using~(b).
\end{proof}

\bigskip

\subsection{Two-sided cells} 
Until now, there is no alternative definition of {\it Calogero-Moser two-sided $c$-cells} 
in terms of Gaudin operators or something related: the only available definition 
is based on the ramification theory of the Calogero-Moser space~\cite[Part~III]{calogero}. 
This depends on the choice of some prime ideal in some Galois closure of some ring extension. 
This choice can be adapted to the choice of the two continuous functions $\g$ and $\xi$ 
and we will follow this choice. 

Moreover, to each Calogero-Moser two-sided $c$-cell $\G$ is associated a subset $\Irr_\G^{c,\calo}(W)$ 
of $\Irr(W)$, which is called a {\it Calogero-Moser $c$-family}. They form a partition of $\Irr(W)$:
\equat\label{eq:two-cm-partition}
\Irr(W) =\dot{\bigcup_\G} ~\Irr_\G^{c,\calo}(W),
\endequat
where $\G$ runs over the set of Calogero-Moser two-sided $c$-cells. 
The following properties are proved 
in~\cite[Theo.~10.2.7,~Prop.~11.3.3,~Prop.~11.4.2]{calogero}:

\bigskip

\begin{prop}\label{prop:two-cm}
Let $\G$ be a Calogero-Moser two-sided $c$-cell and $C$ be a Calogero-Moser left $c$-cell. 
Then:
\begin{itemize}
\itemth{a} $\G$ is a union of Calogero-Moser left 
$c$-cell. Moreover,
$$|\G|=\sum_{\chi \in \Irr_\G^{c,\calo}(W)} \chi(1)^2.$$

\itemth{b} If $C \subset \G$, then every irreducible 
component of $\chi_C^{c,\calo}$ belongs to $\Irr_\G^{c,\calo}(W)$.
\end{itemize}
\end{prop}

\bigskip

It is conjectured~\cite[Conj.~LR]{calogero} that the analogue of Conjecture~\ref{conj:kl=cm} 
also holds for two-sided cells:

\bigskip

\begin{quotation}
\begin{conj}\label{conj:kl=cm-2}
Calogero-Moser two-sided $c$-cells coincide with Kazhdan-Lusztig two-sided $c$-cells. Moreover, 
if $\G$ is one of these, then $\Irr_\G^{c,\calo}(W)=\Irr_\G^{c,\kl}(W)$.
\end{conj}
\end{quotation}

\bigskip

This Conjecture holds in type $A$ (see~\cite{BGW}). In this paper, we 
prove it for dihedral groups whenever $c$ is constant.

\bigskip

\section{Dihedral groups, equal parameters}

\medskip

\boitegrise{{\bf Hypothesis and notation.} {\it From now on, and until the end 
of this paper, we assume that $V=\RM^2$, endowed with its canonical Euclidean 
structure, and we denote by $(e_1,e_2)$ its canonical basis (which is an orthonormal 
basis). We fix a natural number $d \ge 3$ and, if $k \in \ZM$, we set 
$$\a_k=\cos\Bigl(\cfrac{k\pi}{d}\Bigr) e_1 + \sin\Bigl(\cfrac{k\pi}{d}\Bigr) e_2$$
and $s_k=s_{\a_k}$. We assume that $\D=\{\a_0,\a_{d-1}\}$, 
so that $S=\{s_0,s_{d-1}\}$. 
For simplification, we set $s=s_0$ and $s'=s_{d-1}$.}}{0.75\textwidth}

\bigskip

Then $W=\langle s,s' \rangle$ is the dihedral group of order 
$2d$, and $d$ is the order of $ss'$. Moreover, $\a_{k+d}=-\a_k$ and 
$$\Phi=\{\a_k~|~k \in \ZM\}=\{\a_k~|~0 \le k \le 2d-1\}$$
$$\REF(W)=\{s_k~|~k \in \ZM\}=\{s_k~|~0 \le k \le d-1\}.\leqno{\text{and}}$$
Moreover,
$$\Phi^+=\{\a_k~|~0 \le k \le d-1\}.$$
We aim to prove Conjectures~\ref{conj:simple} and~\ref{conj:kl=cm} whenever 
$c$ is constant. The conjecture holds for $c=0$ by Example~\ref{ex:t=0} 
and the remark following Conjecture~\ref{conj:kl=cm}. 
Thus we may assume that $c$ is constant and non-zero and, by Remark~\ref{rem:souplesse}, 
that $c=1$, the constant 
function with value $1$.

\bigskip

\subsection{Elements, characters}
Recall that $s=s_0$ and $s'=s_{d-1}$. 
For $k \ge 0$, we set
$$\s_k=\underbrace{ss's\cdots}_{\text{$k$ terms}}\qquad\text{and}\qquad
\s_k'=\underbrace{s'ss'\cdots}_{\text{$k$ terms}}.$$
Then $\s_0=\s_0'=1$, $\s_d=\s_d'=w_0$ and
\equat\label{eq:elements}
W=\{1,\s_1,\s_1',\s_2,\s_2',\dots,\s_{d-1},\s_{d-1}',w_0\}.
\endequat
Let $W^+=\langle ss' \rangle$. It is a normal cyclic subgroup of order $d$ of $W$. 
We fix a primitive $d$-th root of unity $\o$ and, 
for $k \in \ZM$, we denote by $\th_k : W^+ \to \CM^\times$ the linear character such that 
$\th_k(ss')=\o^k$. We set 
$$\chi_k=\Ind_{W^+}^W \th_k.$$
If $d$ is even, we denote by $\e_s$ (resp. $\e_{s'}$) the linear character of $W$ 
such that $\e_s(s)=-1=-\e_s(s')$ (resp. $\e_{s'}(s)=1=-\e_{s'}(s')$). Then 
$\chi_k=\chi_{-k}=\chi_{k+d}$ and
\equat\label{eq:irr}
\Irr(W)=
\begin{cases}
\{\unb_W,\e,\chi_1,\dots,\chi_{(d-1)/2}\} & \text{if $d$ is odd,}\\
\{\unb_W,\e,\e_s,\e_{s'},\chi_1,\dots,\chi_{(d-2)/2}\} & \text{if $d$ is even.}
\end{cases}
\endequat
We conclude this subsection with a fact that will be useful for our purpose: 
if $k \ge 0$, then
\equat\label{eq:difference}
(\s_k-\s_k') \sum_{s \in \REF(W)} s = 0.
\endequat
\begin{proof}
Note that $\REF(W)=W \setminus W^+$. As it is clear that 
$(\s_k-\s_k')\sum_{w \in W} w =0$, it is sufficient to 
prove that $(\s_k-\s_k')\sum_{w \in W^+} w =0$. But 
$$(\s_k-\s_k')\sum_{w \in W^+} w = \sum_{w \in \s_k W^+} w - \sum_{w \in \s_k' W^+} w,$$
so the result follows from the fact that $\s_kW^+=\s_k'W^+$.
\end{proof}

\bigskip

\subsection{Kazhdan-Lusztig cells} 
Let 
$$C_s=\{\s_1^{-1},\s_2^{-1},\dots,\s_{d-1}^{-1}\}\qquad\text{and}\qquad
C_{s'}=\{\s_1^{\prime -1},\s_2^{\prime -1},\dots,\s_{d-1}^{\prime -1}\}.$$
The Kazhdan-Lusztig left $1$-cells as well as the Kazhdan-Lusztig 
$1$-cellular characters are easily computed (see for instance~\cite[Chap.~21]{livre}):

\bigskip

\begin{prop}\label{prop:kl-diedral}
The Kazhdan-Lusztig left $1$-cells are
$$\{1\},\quad\{w_0\},\quad C_s \quad\text{and}
\quad C_{s'}.$$
Moreover:
\begin{itemize}
\itemth{a} If $d$ is odd, then 
$$\chi_{C_s}^{1,\kl}=\chi_{C_{s'}}^{1,\kl}=\chi_1+\cdots + \chi_{(d-1)/2}.$$

\itemth{b} 
If $d$ is even, then
$$
\chi_{C_s}^{1,\kl}=\e_{s'} + \chi_1+\cdots + \chi_{(d-2)/2} \quad\text{and}\quad  
\chi_{C_{s'}}^{1,\kl}=\e_s + \chi_1 + \cdots + \chi_{(d-2)/2}.$$
\end{itemize}
\end{prop}

\bigskip

\subsection{Calogero-Moser cells}
The main result of our paper is the following:

\bigskip

\begin{theo}\label{theo:main}
Conjectures~\ref{conj:simple-facile},~\ref{conj:kl=cm} and~\ref{conj:kl=cm-2} hold whenever 
$W$ is dihedral and $c$ is constant.
\end{theo}

\bigskip

The rest of this section is devoted to the proof of this Theorem.

\bigskip

\subsection{Preliminaries}
We use the flexibility of the definition of left cells explained in Remark~\ref{rem:souplesse}. 
We take 
$$v_1=v_2=\sin\Bigl(\cfrac{\pi}{2d}\Bigr) e_1 + \cos\Bigl(\cfrac{\pi}{2d}\Bigr) e_2,
\qquad\g(t)=t\qquad\text{and}\qquad \xi(t)=1-t$$
(one can check that $v_0$ is a positive multiple of $v_1$). We also set 
$$v_1^\perp=\cos\Bigl(\cfrac{\pi}{2d}\Bigr) e_1 - \sin\Bigl(\cfrac{\pi}{2d}\Bigr) e_2,$$
so that $(v_1,v_1^\perp)$ is an orthonormal basis of $V$. We have 
\equat\label{eq:scalaire}
\begin{cases}
(v_1,\a_k)=\DS{\sin\Bigl(\cfrac{(1+2k)\pi}{2d}\Bigr)},\\
(v_1^\perp,\a_k)=\DS{\cos\Bigl(\cfrac{(1+2k)\pi}{2d}\Bigr)}\\
\end{cases}
\endequat
for all $k \ge 0$. Also
\equat\label{eq:action sk}
\begin{cases}
\s_k(v_1)=\DS{\sin\Bigl(\cfrac{(1-2k)\pi}{2d}\Bigr) e_1 + \cos\Bigl(\cfrac{(1-2k)\pi}{2d}\Bigr) e_2,} \\
\\
\DS{\s_k'(v_1)=\sin\Bigl(\cfrac{(1+2k)\pi}{2d}\Bigr) e_1 + \cos\Bigl(\cfrac{(1+2k)\pi}{2d}\Bigr) e_2.}
\end{cases}
\endequat

\medskip

\begin{proof}[Proof of~\eqref{eq:action sk}]
This is easily checked by induction on $k$ using the fact that $\s_k=s\s_{k-1}'$ and $\s_k'=s'\s_{k-1}$ 
for $k \ge 1$.
\end{proof}

\medskip

An immediate consequence is that
\equat\label{eq:sk scalaire}
\begin{cases}
(\s_k(v_1),v_1)=(\s_k'(v_1),v_1)=\DS{\cos\Bigl(\cfrac{k\pi}{d}\Bigr)},\\
\\
(\s_k(v_1),v_1^\perp)=\DS{-\sin\Bigl(\cfrac{k\pi}{d}\Bigr)},\\
\\
(\s_k'(v_1),v_1^\perp)=\DS{\sin\Bigl(\cfrac{k\pi}{d}\Bigr)}.\\
\end{cases}
\endequat

\bigskip

Now, we set $\r=ss'$ (so that $\r$ is the rotation with angle $2\pi/d$). Recall that 
$\r$ generates $W^+$. Also
\equat\label{eq:rho}
s_ks_l=\r^{l-k}
\endequat
for all $k$, $l \in \ZM$. 

\bigskip

\subsection{Proof of Conjecture~\ref{conj:simple-facile}}\label{sub:simple}
Now, we set 
$$A=D_{-v_1}^{0,v_1,v_1},\quad B=D_{-v_1}^{1,v_1,0},\quad
A^\perp=D_{v_1^\perp}^{0,v_1,v_1}\quad\text{and}\quad B^\perp=D_{v_1^\perp}^{1,v_1,0}.$$
In particular,
$$D_{-v_1}^{b \cdot 1,v_1,av_1}=aA+bB\qquad\text{and}\qquad 
D_{v_1^\perp}^{b \cdot 1,v_1,av_1}=aA^\perp+bB^\perp.$$
Now, if $1 \le k \le d-1$, then
\equat\label{eq:propres}
\begin{cases}
\text{$\s_k^{-1}-\s_k^{\prime -1}$ is an eigenvector of $aA+bB$ for the eigenvalue $-a\cos(k\pi/d)$,}\\
\text{$\s_k^{-1}-\s_k^{\prime -1}$ is not an eigenvector of $aA^\perp+bB^\perp$ if $a > 0$.}
\end{cases}
\endequat

\medskip

\begin{proof}[Proof of~\eqref{eq:propres}]
Note that $\s_k^{-1}-\s_k^{\prime -1}=\pm(\s_k-\s_k')$. 
By~\eqref{eq:sk scalaire}, we have $A(\s_k^{-1}-\s_k^{\prime -1})=-\cos(k\pi/d)(\s_k^{-1}-\s_k^{\prime -1})$. 
Moreover, $B(w)=w \sum_{s \in \REF(W)}s$, so $A(\s_k^{-1}-\s_k^{\prime -1})=0$ by~\eqref{eq:difference}. 
The first assertion follows.

Still by~\eqref{eq:sk scalaire}, 
the coefficient of $\s_k^{-1}$ (or $\s_k^{\prime -1}$) in 
$(aA^\perp+bB^\perp)(\s_k^{-1}-\s_k^{\prime -1})$ 
is equal to $-a\sin(k\pi/d)$, so this proves the second statement because $\sin(k\pi/d) \neq 0$.
\end{proof}

\medskip

After these preliminaries, we are ready to prove the theorem.
So let us first prove that Conjecture~\ref{conj:simple-facile} holds. 
Let us assume that $a > 0$. Let $\r_{a,b}$ denote the largest eigenvalue of 
$aA+bB$ (as in the proof of Proposition~\ref{prop:calo-w0}(c)). Then it has 
multiplicity $1$. By the proof of Proposition~\ref{prop:calo-w0}(b), 
$-\r_{a,b}$ is the smallest eigenvalue of $aA+bB$, and it has multiplicity $1$. 
We denote by $E_\maxi(a,b)$ (resp. $E_\mini(a,b)$) the $\r_{a,b}$-eigenspace 
(resp. the $-\r_{a,b}$-eigenspace) of $aA+bB$. 

On the other hand, if $1 \le k \le d-1$, 
it follows from~\eqref{eq:propres} that the vector space generated by $\s_k^{-1}-\s_k^{\prime -1}$ 
and its image by $aA^\perp+bB^\perp$ is contained in the $-a\cos(k\pi/d)$-eigenspace 
of $aA+bB$. So this eigenspace (let us denote it by $E_k(a,b)$) has dimension $\ge 2$. 
Since $a \neq 0$, $a \cos(k\pi/d) \neq a\cos(l\pi/d)$ if $1 \le k < l \le d-1$. Therefore, 
the vector space 
$$E_\maxi(a,b) \oplus E_\mini(a,b) \oplus E_1(a,b) \oplus E_2(a,b) \oplus \cdots \oplus E_{d-1}(a,b)$$
has dimension $\ge 2 + 2(d-1)=2d = \dim \RM W$. So this proves that
\equat\label{eq:eab}
\RM W=E_\maxi(a,b) \oplus E_\mini(a,b) \oplus E_1(a,b) \oplus E_2(a,b) \oplus \cdots \oplus E_{d-1}(a,b)
\endequat
and that 
\equat\label{eq:sim 2}
\dim E_k(a,b)=2
\endequat
for $1 \le k \le d-1$. 
This describes the diagonalization of $aA+bB$. But now, the second statement 
of~\eqref{eq:propres} shows that $aA^\perp + bB^\perp$ has two distinct eigenvalues on each 
$E_k(a,b)$. So the family $\DC^{a \cdot 1,v_1,bv_1}$ has simple spectrum as soon as $a > 0$. 
This is exactly Conjecture~\ref{conj:simple-facile}.

\medskip

\subsection{Cellular characters}
Let us define two elements of the group algebra $\RM W$ by 
$$\bG=\sum_{k=0}^{d-1} s_k\qquad\text{and}\qquad
\bG^\perp=\sum_{k=0}^{d-1} \cot\Bigl(\cfrac{(2k+1)\pi}{2d}\Bigr) s_k.$$
Then $B$ (resp. $B^\perp$) is the right multiplication by $\bG$ (resp. $\bG^\perp$). 
For $\chi \in \Irr(W)$, we denote by $e_\chi \in \Zrm(\RM W)$ 
the corresponding central idempotent
$$e_\chi=\cfrac{\chi(1)}{|G|} \sum_{w \in W} \chi(w^{-1}) w.$$
Then 
$$\RM W = \bigoplus_{\chi \in \Irr(W)} \RM W e_\chi$$
and $\bG=d(e_{\unb_W}-e_\e)$. We denote by $E_d$ (resp. $E_{-d}$) the $d$-eigenspace 
(resp. $(-d)$-eigenspace) of $B^\perp$. Then:

\bigskip

\begin{lem}\label{lem:diagonalisation}
With the above notation, we have:
\begin{itemize}
 \itemth{a} $\RM W = \RM e_{\unb_W} \oplus \RM e_\e \oplus E_d \oplus E_{-d}$ and
 $$E_d \oplus E_{-d}=\bigoplus_{\chi \in \Irr(W)\setminus \{\unb_W,\e\}} \RM W e_\chi.$$
 
 \itemth{b} $B$ acts on $\RM e_{\unb_W}$ (resp. $\RM e_\e$, resp. $E_d$, resp. $E_{-d}$) 
 by multiplication by $d$ (resp. $-d$, resp. $0$, resp. $0$).
 
 \itemth{c} $B^\perp$ acts on $\RM e_{\unb_W}$ and $\RM e_\e$ by multiplication by $0$.
\end{itemize}
\end{lem}

\bigskip

\begin{proof}
First, it follows from~\eqref{eq:rho} and~\eqref{eq:cot} that 
$$(\bG^\perp)^2 = d^2 - d \sum_{w \in W^+} w.$$
Moreover, 
$$\bG^2=d \sum_{w \in W^+} w.$$
But
$$e_\chi \Bigl(\sum_{w \in W^+} w\Bigr) = 
\begin{cases}
d e_\chi & \text{if $\chi \in \{\unb_W,\e\}$,}\\
0 & \text{if $\chi \not\in \{\unb_W,\e\}$.}
\end{cases}$$
This proves the lemma because $B$ and $B^\perp$ are diagonalizable.
\end{proof}

\bigskip

This shows that 
$$\spectrum^{1,v_1,0} =\{dv_1^*,-dv_1^*,d(v_1^\perp)^*,-d(v_1^\perp)^*\}$$
and that the decomposition in~(a) is the decomposition~\eqref{eq:decomposition} 
for $(c,v,v')=(1,v_1,0)$. Let us now determine the corresponding cellular 
characters:

\bigskip

\begin{lem}\label{lem:cellular diedral}
If $d$ is odd, then $E_d$ and $E_{-d}$ both afford the character 
$\chi_1+\chi_2+\cdots +\chi_{(d-1)/2}$. 

If $d$ is even, then the $\RM W$-module $E_d$ (resp. $E_{-d}$) affords the character 
$\e_{s'}+\chi_1+\chi_2+\cdots +\chi_{(d-2)/2}$ (resp. $\e_s+\chi_1+\chi_2+\cdots +\chi_{(d-1)/2}$). 
\end{lem}

\bigskip

\begin{proof}
Let $\eta$ denote the automorphism of $W$ exchanging $s$ and $s'$. Then 
$\lexp{\eta}{\bG^\perp}=-\bG^\perp$. So $\eta$ exchanges $E_d$ and $E_{-d}$. 
Moreover, $\lexp{\eta}{\chi_j}=\chi_j$ for any $j$. 

But, if $d$ is odd, then 
$E_d \oplus E_{-d}$ affords the character $2(\chi_1+\chi_2+\cdots +\chi_{(d-1)/2})$. 
So this forces that both $E_d$ and $E_{-d}$ afford the character 
$\chi_1+\chi_2+\cdots +\chi_{(d-1)/2}$. 

Now, if $d$ is even, then $E_d \oplus E_{-d}$ affords the character 
$\e_s+\e_{s'}+2(\chi_1+\chi_2+\cdots +\chi_{(d-1)/2})$. So this forces 
that $E_d$ affords the character $\e_s+\chi_1+\chi_2+\cdots +\chi_{(d-2)/2})$ 
or $\e_{s'}+\chi_1+\chi_2+\cdots +\chi_{(d-2)/2})$ (and $E_{-d}$ affords the other). 
To determine which one it is, we just need to compute $e_{\e_s} \bG^\perp$. 
But $e_{\e_s} \bG^\perp=\l e_{\e_s}$ where 
$$\l=\sum_{k=0}^{d-1} \cot\Bigl(\cfrac{(2k+1)\pi}{2d}\Bigr) \e_s(s_k).$$
Since $\e_s(s_k)=(-1)^{k+1}$, we get from~\eqref{eq:cot d} that $\l=-d$. 
This proves the result.
\end{proof}

\bigskip

Lemmas~\ref{lem:diagonalisation} and~\ref{lem:cellular diedral}, together 
with Proposition~\ref{prop:kl-diedral}, prove that 
the list of Calogero-Moser $1$-cellular characters equals the list 
of Kazhdan-Lusztig $1$-cellular characters. 

\bigskip

\subsection{Left cells}
Assume that $a \neq 0$. 
If $1 \le k \le d-1$, then $E_k(a,b)$ has dimension $2$ (see~\S\ref{sub:simple}). 
Let us determine the action of $aA^\perp + bB^\perp$ on $E_k(a,b)$. 
Let $r_k(a,b)=(aA^\perp + bB^\perp)(\s_k^{-1}-\s_k^{\prime -1}) \in E_k(a,b)$. Then
$$E_k(a,b) = \RM (\s_k^{-1}-\s_k^{\prime -1}) \oplus \RM r_k(a,b).$$
So this means that we need to determine $(aA^\perp + bB^\perp)(r_k(a,b))$, 
i.e. we need to determine the two real numbers $z$, $z'$ such that 
$$(aA^\perp + bB^\perp)(r_k(a,b))=z (\s_k^{-1}-\s_k^{\prime -1}) + z' r_k(a,b).$$
For this, we will work modulo $F_k=\bigoplus_{w \in W \setminus \{\s_k,\s_k'\}} \RM w$.
Then
$$r_k(a,b) \equiv -a \sin\Bigl(\cfrac{k\pi}{d}\Bigr)(\s_k^{-1}+\s_k^{\prime -1}) \mod F_k.$$
On the other hand,
$$(aA^\perp + bB^\perp)^2=a^2 (A^\perp)^2 + b^2 (B^\perp)^2 + ab (A^\perp B^\perp + B^\perp A^\perp).$$
So
$$(aA^\perp + bB^\perp)(r_k(a,b)) \equiv 
a^2\sin^2\Bigl(\cfrac{k\pi}{d}\Bigr)(\s_k^{-1}-\s_k^{\prime -1})
+ b^2 (B^\perp)^2(\s_k^{-1}-\s_k^{\prime -1}) \mod F_k.$$
But the formula for $(\bG^\perp)^2$ given in the proof of Lemma~\ref{lem:diagonalisation} 
shows that 
$$(B^\perp)^2(\s_k^{-1}-\s_k^{\prime -1})=d^2(\s_k^{-1}-\s_k^{\prime -1}).$$
All this together shows that 
$$z=a^2\sin^2\Bigl(\cfrac{k\pi}{d}\Bigr)+b^2d^2\qquad\text{and}\qquad z'=0.$$
So we have proved the following fact:

\bigskip

\begin{lem}\label{lem:vp}
Assume that $a \neq 0$. Then 
the restriction of $aA^\perp + bB^\perp$ to $E_k(a,b)$ has two eigenvalues, 
namely 
$$\pm \sqrt{a^2\sin^2\Bigl(\cfrac{k\pi}{d}\Bigr)+b^2d^2}.$$
\end{lem}

\bigskip

Coming back to the family $\DC^{t\cdot c,v_1,(1-t)v_1}$ (which corresponds 
to the case $a=1-t$ and $b=t$), this shows that 
\equat\label{eq:vp-diedral}
\begin{cases}
\l_{\s_k^{-1}}(t)=
\DS{-(1-t)\cos\Bigl(\cfrac{k\pi}{d}\Bigr)v_1^* - 
\sqrt{(1-t)^2\sin^2\Bigl(\cfrac{k\pi}{d}\Bigr)+t^2d^2}}~ (v_1^\perp)^*,\\
~\\
\l_{\s_k^{\prime -1}}(t)=
\DS{-(1-t)\cos\Bigl(\cfrac{k\pi}{d}\Bigr)v_1^* + \sqrt{(1-t)^2\sin^2\Bigl(\cfrac{k\pi}{d}\Bigr)+t^2d^2}} 
~(v_1^\perp)^*,\\
\end{cases}
\endequat
for $1 \le k \le d-1$. Taking the limit at $t=1$, and using the cases $1$ and $w_0$ 
solved in Proposition~\ref{prop:calo-w0}, we get (also thanks to Lemma~\ref{lem:cellular diedral}):

\bigskip

\begin{prop}\label{prop:cm-diedral}
The Calogero-Moser left $1$-cells are
$$\{1\},\quad\{w_0\},\quad C_s \quad\text{and}
\quad C_{s'}.$$
Moreover:
\begin{itemize}
\itemth{a} If $d$ is odd, then 
$$\chi_{C_s}^{1,\calo}=\chi_{C_{s'}}^{1,\calo}=\chi_1+\cdots + \chi_{(d-1)/2}.$$

\itemth{b} 
If $d$ is even, then
$$
\chi_{C_s}^{1,\calo}=\e_{s'} + \chi_1+\cdots + \chi_{(d-2)/2} \quad\text{and}\quad  
\chi_{C_{s'}}^{1,\calo}=\e_s + \chi_1 + \cdots + \chi_{(d-2)/2}.$$
\end{itemize}
\end{prop}

\bigskip

The proof of Conjecture~\ref{conj:kl=cm-2} in this case is complete.

\bigskip

\subsection{Two-sided cells} 
The Kazhdan-Lusztig two-sided $1$-cells are 
$$\{1\}, \quad \G=C_s \cup C_{s'} \quad \text{and}\quad \{w_0\}.$$
The Kazhdan-Lusztig $1$-families are given by
$$\Irr_{\{1\}}^{1,\kl}(W)=\{\e\}
,\quad \Irr_\G^{1,\kl}(W)=\Irr(W) \setminus \{\unb_W,\e\}
\quad\text{and}\quad \Irr_{\{w_0\}}^{1,\kl}(W)=\{\unb_W\}.$$
The Calogero-Moser $1$-families have been computed by Bellamy~\cite{bellamy these} 
(see also~\cite[Table~5.1]{bonnafe diedral}) and coincide with the 
Kazhdan-Lusztig $1$-families. Therefore, the fact that Conjecture~\ref{conj:kl=cm-2} 
holds in this case follows immediately from Proposition~\ref{prop:two-cm}. 

The proof of Theorem~\ref{theo:main} is complete.

\section{Complements}

\medskip

\subsection{Other eigenvalues}
Keep the notation of the proof of Theorem~\ref{theo:main}. 
In this proof, it was unnecessary 
to determine the explicit value of the largest eigenvalue $\r_{a,b}$ of $aA+bB$ and the 
action of $aA^\perp+bB^\perp$ on $E_\mini(a,b)$ and $E_\maxi(a,b)$. This can easily be done, 
as shown in this subsection. 

First, note that $\Tr(AB)=0$ because the diagonal coefficients are equal to zero. 
Also, 
$$
\Tr(A^2)=\sum_{k=0}^{2d-1} \cos^2\Bigl(\frac{k\pi}{d}\Bigr)
\qquad\text{and}\qquad 
\Tr(B^2)=2d^2.
$$
The last equality comes from the fact that the coefficient of $1$ in $(\sum_{s \in \REF(W)} s)^2$ 
is equal to $d$. Therefore, 
$$
\Tr((aA+bB)^2)=2a^2 \sum_{k=0}^{d-1} \cos^2\Bigl(\frac{k\pi}{d}\Bigr)+ 2b^2 d^2
$$
Since the eigenvalue of $aA+bB$ on $E_\mini(a,b)$ is $-\r_{a,b}$, 
the decomposition~\eqref{eq:eab} implies that 
$$\Tr((aA+bB)^2)=2 \r_{a,b}^2 + 2a^2 \sum_{k=1}^{d-1} \cos^2\Bigl(\frac{k\pi}{d}\Bigr).$$
This shows that
\equat\label{eq:rhoab}
\r_{a,b}=\sqrt{a^2+d^2b^2}
\endequat

It remains to determine the action of $aA^\perp+bB^\perp$ on $E_\mini(a,b)$ and $E_\maxi(a,b)$. 
Both spaces have dimension $1$, so this action is by a scalar: if $\r$ denotes the scalar 
by which $aA^\perp+bB^\perp$ acts on $E_\maxi(a,b)$, then $-\r$ is the scalar 
by which $aA^\perp+bB^\perp$ acts on $E_\mini(a,b)$ (again by the proof of Proposition~\ref{prop:calo-w0}). 
Now, note that $\Tr(A^\perp B^\perp)=0$ because the diagonal coefficients are equal to zero. 
Also, 
$$
\Tr((A^\perp)^2)=\sum_{k=0}^{2d-1} \sin^2\Bigl(\frac{k\pi}{d}\Bigr)
\qquad\text{and}\qquad 
\Tr((B^\perp)^2)=2d^2(d-1).
$$
Therefore, 
$$
\Tr((aA^\perp+bB^\perp)^2)=2a^2 \sum_{k=1}^{d-1} \sin^2\Bigl(\frac{k\pi}{d}\Bigr)+ 2b^2 d^2(d-1).
$$
The last equality follows from the formula for $(\bG^\perp)^2$ given in the proof 
of Lemma~\ref{lem:diagonalisation}. But the decomposition~\eqref{eq:eab} implies that 
$$\Tr((aA^\perp+bB^\perp)^2)=
2\r^2 + 2 \sum_{k=1}^{d-1} \Bigl(a^2\sin^2\Bigl(\frac{k\pi}{d}\Bigr)+b^2d^2\Bigr).$$
This shows that $\r^2=0$ and so
\equat\label{eq:rho-0}
\r=0.
\endequat
Keeping the notation of the proof of Theorem~\ref{theo:main}, this gives 
the following formula for the paths $\l_1(t)$ and $\l_{w_0}(t)$:
\equat\label{eq:lambda-1}
\l_1(t)=\sqrt{(1-t)^2+d^2t^2}~v_1^*
\qquad\text{and}\qquad \l_{w_0}(t)=-\sqrt{(1-t)^2+d^2t^2}~v_1^*.
\endequat

\subsection{Some pictures}
\begin{figure}
\begin{centerline}{
\begin{tabular}{ccc}
\includegraphics[scale=0.4]{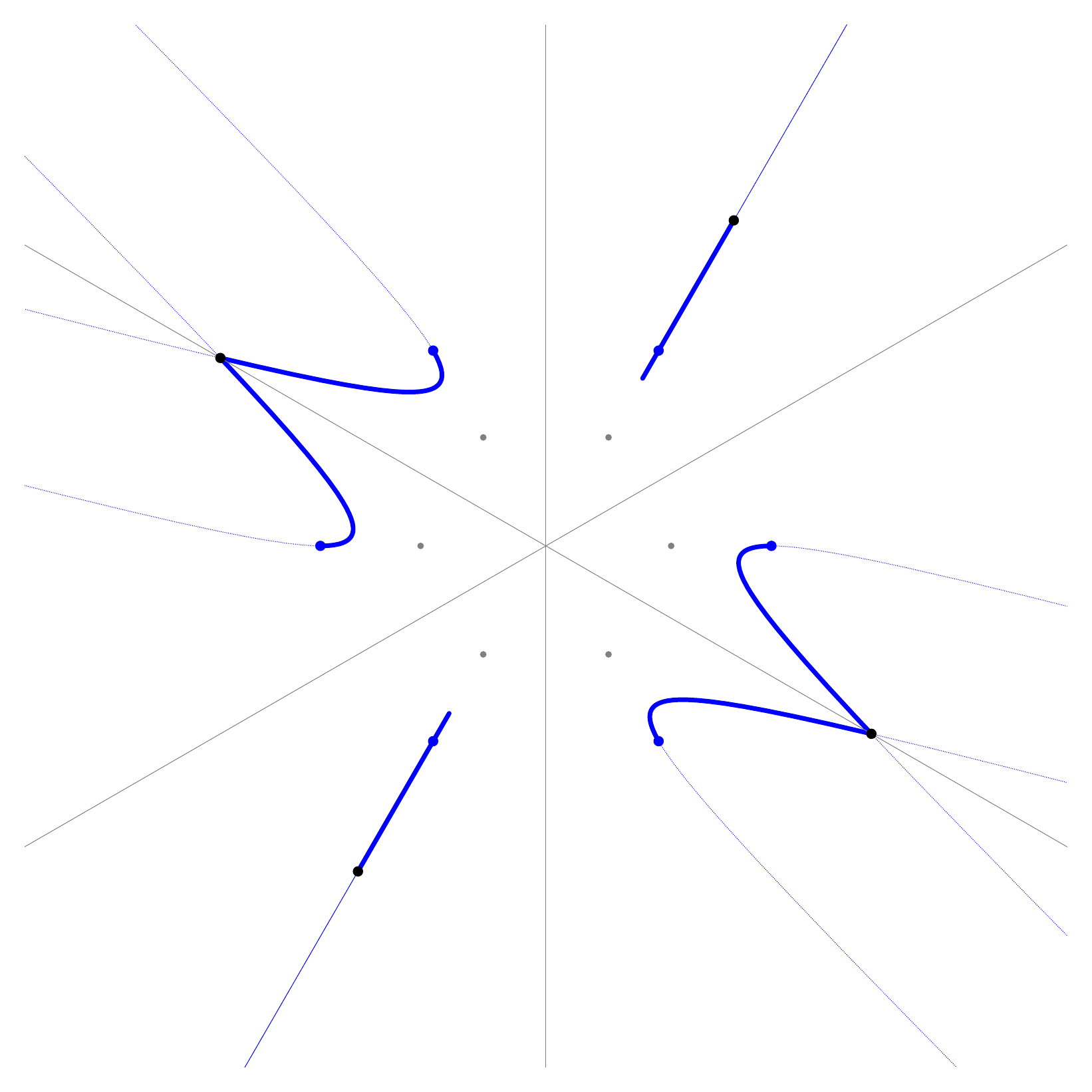}&~&
\includegraphics[scale=0.4]{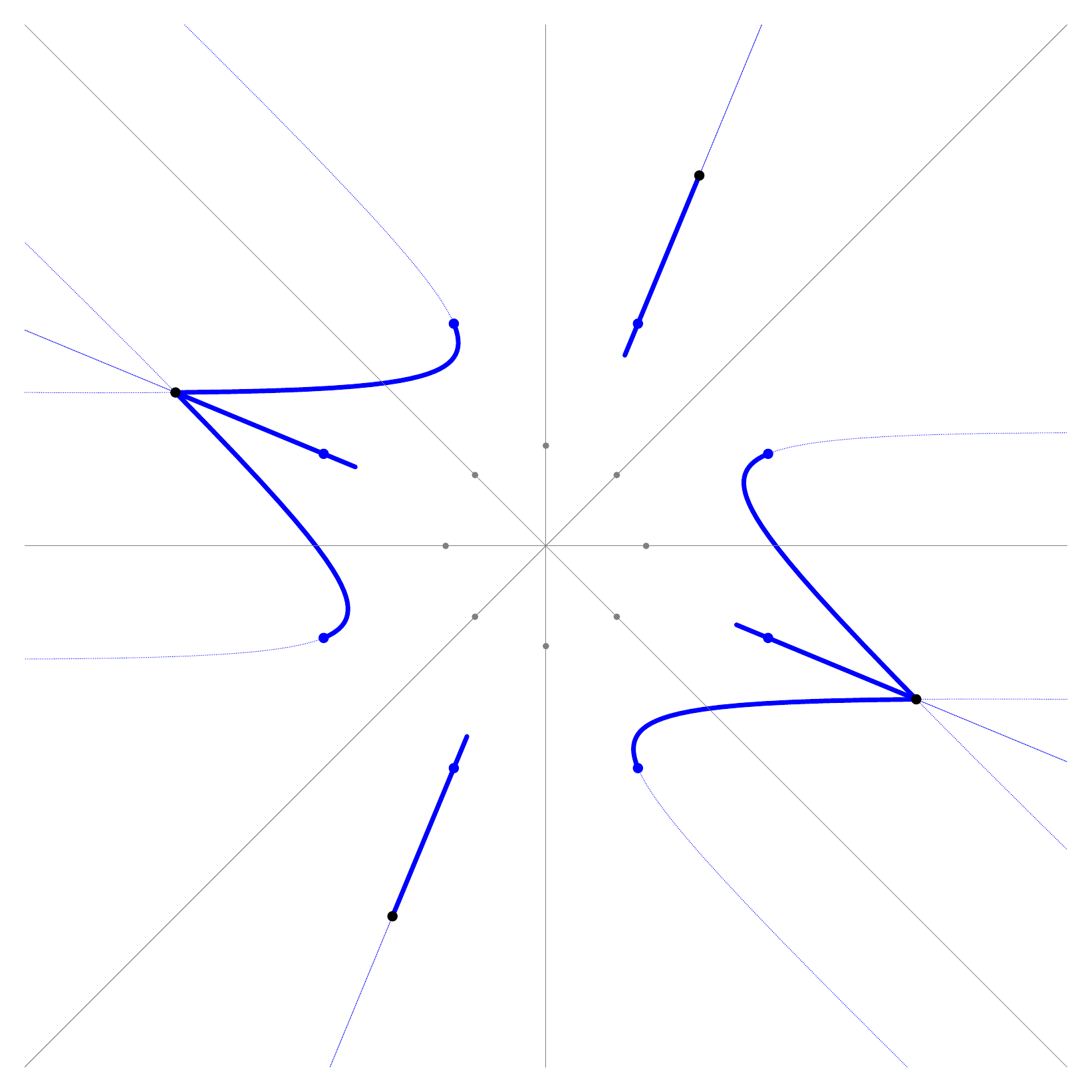}\\
$d=3$ && $d=4$ \\
&&\\
\includegraphics[scale=0.4]{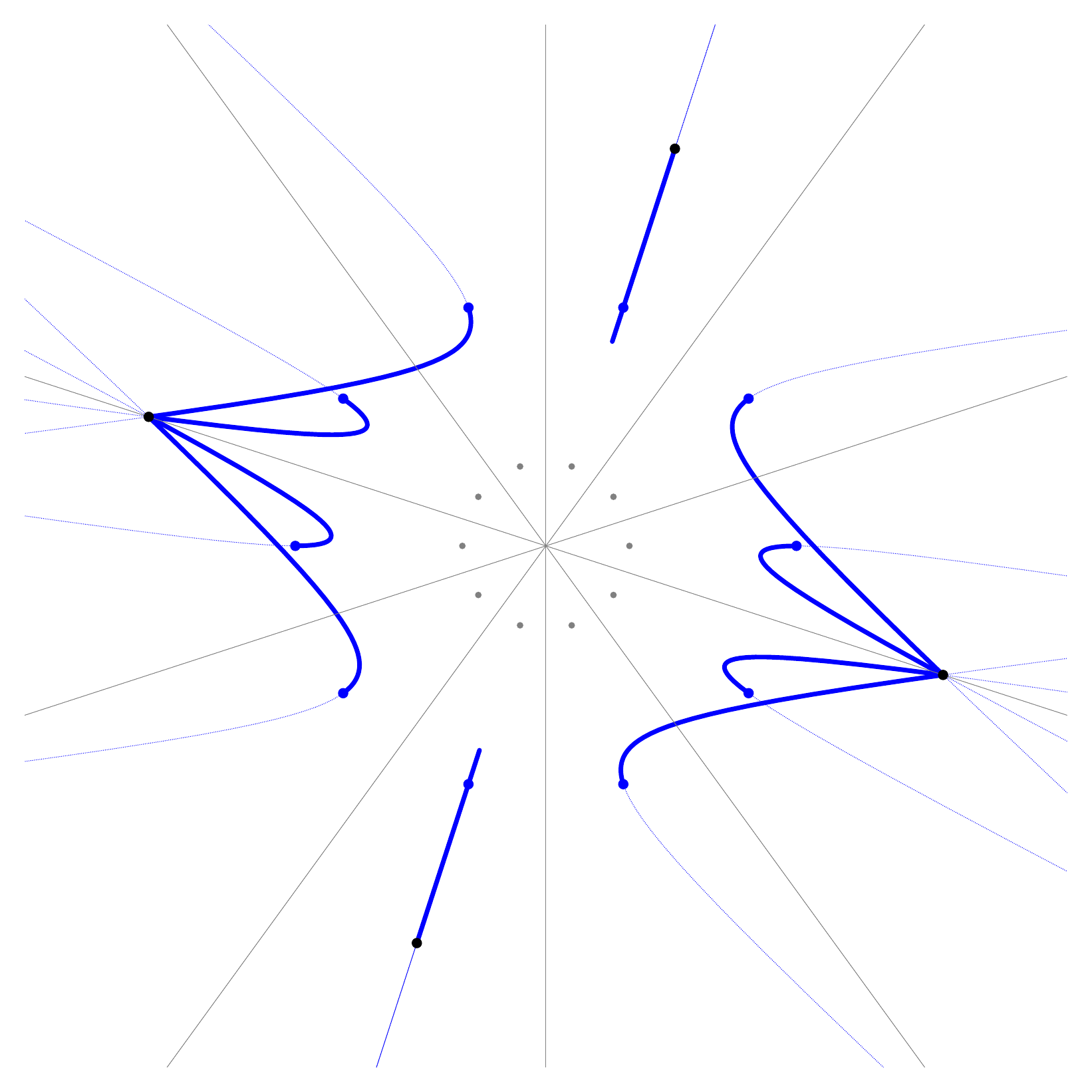}&&
\includegraphics[scale=0.4]{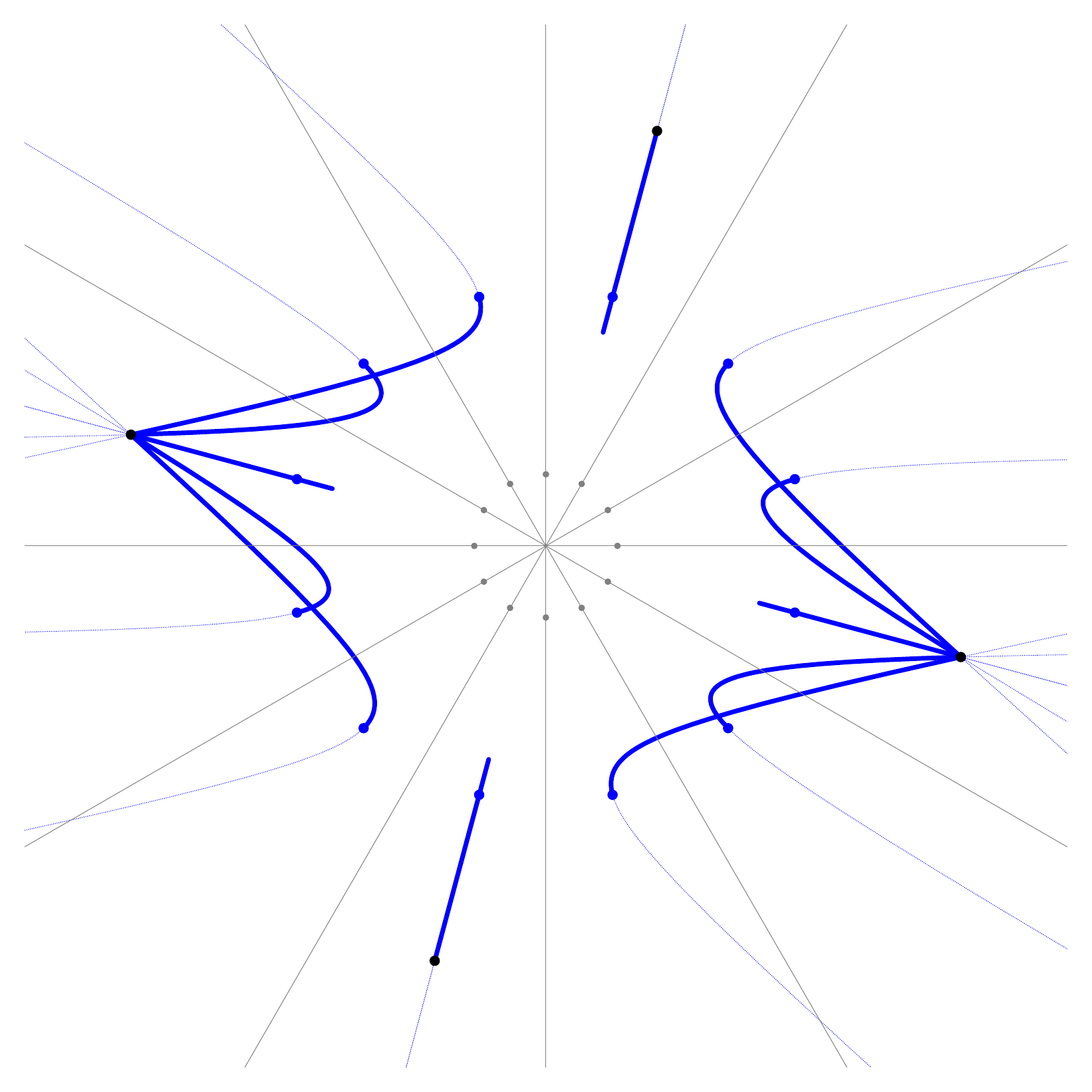}\\
$d=5$ && $d=6$ \\
&&\\
\includegraphics[scale=0.4]{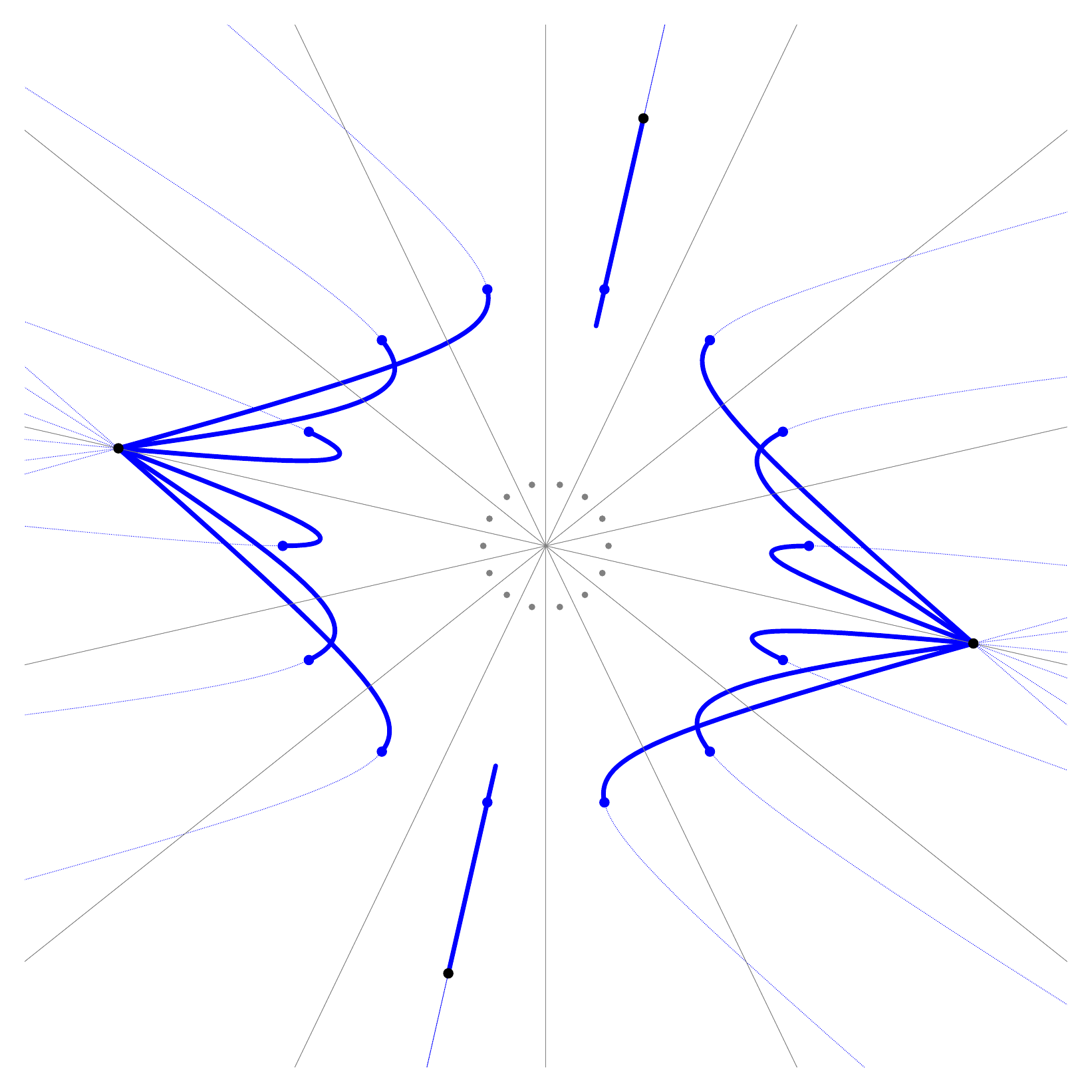}&&
\includegraphics[scale=0.4]{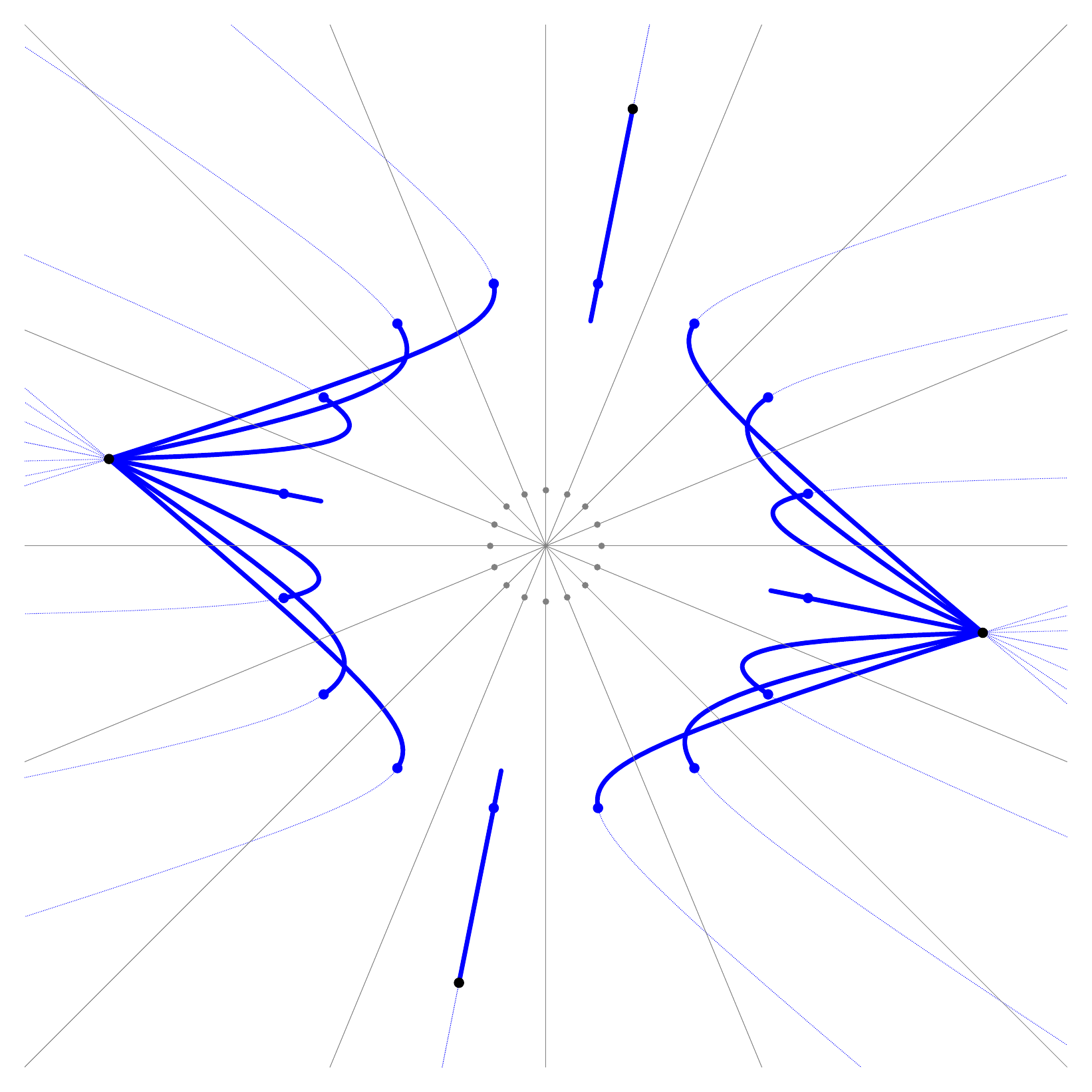}\\
$d=7$ && $d=8$ \\
\end{tabular}}
\end{centerline}
\caption{Paths $(\l_w)_{w \in W}$ for $3 \le d \le 8$}\label{fig:paths}
\end{figure}

We provide in Figure~\ref{fig:paths} 
some pictures of the paths described in~\eqref{eq:vp-diedral} and~\eqref{eq:lambda-1}, 
whenever $d \in \{3,4,5,6,7,8\}$. In these pictures, we have identified $V$ and $V^*$ 
as all along the paper. The gray points represent the roots, the gray lines represent 
the reflecting hyperplanes, the blue dots are the points $w(3dv_1/5)$ (the reason 
for choosing $3dv_1/5$ is to have a better view of what happens), the black dots 
represent the spectrum of the family $\DC^{1,v_1,0}$ (i.e. they are in bijection 
with Calogero-Moser left cells), the blue thick curves are the 
paths $\l_w(t)$ (renormalized as above, i.e. for the family $\DC^{1,tv_1,3d(1-t)v_1/5}$) 
for $t \in [0,1]$ and the blue thin lines 
are the extensions of these paths for arbitrary values of $t$. 

\setcounter{section}{0}
\renewcommand\thesection{\Roman{section}}
\def\sectionname{Appendix}

\section{Trigonometric identities}

\medskip

We aim to prove that
\equat\label{eq:cot}
\sum_{k=0}^{d-1} \cot\Bigl(\cfrac{(2k+1)\pi}{2d}\Bigr)\cot\Bigl(\cfrac{(2k+2j+1)\pi}{2d}\Bigr)=
\begin{cases}
d^2-d & \text{if $j=0$,}\\
-d & \text{if $1 \le j \le d-1$.}
\end{cases}
\endequat

\bigskip

\begin{proof}
For $n \ge 1$, let 
$f(n)=\sum_{k=1}^{n-1} \cot^2\Bigl(\cfrac{k\pi}{n}\Bigr)$. The fact that 
$f(n)=\cfrac{(n-1)(n-2)}{3}$ goes back (at least) to Cauchy. Now,
$$\sum_{k=0}^{d-1} \cot^2\Bigl(\cfrac{(2k+1)\pi}{2d}\Bigr)=f(2d)-f(d),$$
so the first equality follows easily. 

Assume now that $1 \le j \le d-1$. Since $\cot(x)\cot(y)=\cot(y-x)(\cot(x)-\cot(y))-1$ 
whenever $y \not\equiv x \!\!\mod \pi$, we get
\begin{multline*}
\sum_{k=0}^{d-1} \cot\Bigl(\cfrac{(2k+1)\pi}{2d}\Bigr)\cot\Bigl(\cfrac{(2k+2j+1)\pi}{2d}\Bigr)=\\ 
-d + \cot\Bigl(\cfrac{j\pi}{d}\Bigr) \sum_{k=0}^{d-1} \Bigl(\cot\Bigl(\cfrac{(2k+1)\pi}{2d}\Bigr)
-\cot\Bigl(\cfrac{(2k+2j+1)\pi}{2d}\Bigr)\Bigr).
\end{multline*}
But the sequences $(k+j)_{0 \leqslant k \leqslant d-1}$ and $(k)_{0 \leqslant k \leqslant d-1}$ both 
cover all the integers modulo $d$. 
So the sum of the terms in the last summand vanishes, as desired.
\end{proof}

\bigskip

Note also the following trivial identity:
\equat\label{eq:cot nul}
\sum_{k=0}^{d-1} \cot\Bigl(\cfrac{(2k+1)\pi}{2d}\Bigr)=0
\endequat
(the terms indexed by $k$ and $d-1-k$ are opposite to each other). Also, if $d$ is even, then
\equat\label{eq:cot d}
\sum_{k=0}^{d-1} (-1)^k\cot\Bigl(\cfrac{(2k+1)\pi}{2d}\Bigr)=d.
\endequat

\begin{proof}
Let $\xi=\exp(i\pi/(2d))$ and $\z=\xi^2$. Then $\z^d=-1$ and 
$$\cot\Bigl(\cfrac{(2k+1)\pi}{2d}\Bigr)=\xi^d \cfrac{\x^{2k+1}+\xi^{-2k-1}}{\x^{2k+1}-\xi^{-2k-1}}=
\xi^d \cfrac{\z^{2k+1}+1}{\z^{2k+1}-1}.$$
Since $d$ is even, we can write $d=2m$. Then $\xi^d=\z^m$ and
$$(-1)^k=\z^{2mk}=\z^{m(2k+1)-m}\qquad\text{and}\qquad 
(-1)^k\z^{2k+1}=\z^{(m+1)(2k+1)-m}.$$
Therefore,
$$(-1)^k\cot\Bigl(\cfrac{(2k+1)\pi}{2d}\Bigr)=\cfrac{\z^{m(2k+1)}+\z^{(m+1)(2k+1)}}{\z^{2k+1}-1}.$$
The result follows from~\cite[$(1.10)$]{bonnafe diedral}, specialized at $X=Y=1$.
\end{proof}

\end{document}